\documentclass[12pt]{article}
\usepackage{amsmath,amsfonts,amssymb}

\newtheorem{thm}{Theorem}[section]
\newtheorem{cor}[thm]{Corollary}
\newtheorem{prop}[thm]{Proposition}
\newtheorem{lem}[thm]{Lemma}

\newtheorem{con}[thm]{Conjecture}

\def\qed{\nopagebreak\hfill{\rule{4pt}{7pt}}}
\def\pf{\noindent {\it Proof.} }

\begin{document}

\begin{center}
{\large \bf The srank Conjecture on Schur's $Q$-Functions}
\end{center}

\begin{center}
William Y. C. Chen$^{1}$, Donna Q. J. Dou$^2$, \\
Robert L. Tang$^3$
and
Arthur L. B. Yang$^{4}$\\[6pt]
Center for Combinatorics, LPMC-TJKLC\\
Nankai University, Tianjin 300071, P. R. China\\[5pt]
$^{1}${\tt chen@nankai.edu.cn}, $^{2}${\tt
qjdou@cfc.nankai.edu.cn}, $^{3}${\tt tangling@cfc.nankai.edu.cn},
$^{4}${\tt yang@nankai.edu.cn}
\end{center}


\noindent \textbf{Abstract.} We show that the shifted rank, or
srank, of any partition $\lambda$ with distinct parts equals the
lowest degree of the terms  appearing in the expansion of  Schur's
$Q_{\lambda}$ function in terms of power sum symmetric functions.
This gives an affirmative answer to a conjecture of Clifford. As
pointed out by Clifford, the notion of the srank can be naturally
extended to a skew partition $\lambda/\mu$ as the minimum number of
bars among the corresponding skew bar tableaux. While the srank
conjecture is not valid for skew partitions, we give an algorithm to
compute the srank.


\noindent \textbf{MSC2000 Subject Classification:} 05E05, 20C25

\section{Introduction}

The main objective of this paper is to answer two open problems
raised by Clifford \cite{cliff2005} on sranks of partitions with
distinct parts, skew partitions and Schur's $Q$-functions. For any
partition $\lambda$ with distinct parts, we give a proof of
Clifford's srank conjecture that the lowest degree of the terms in
the power sum expansion of Schur's $Q$-function $Q_{\lambda}$ is
equal to the number of bars in a minimal bar tableaux of shape
$\lambda$. Clifford \cite{cliff2003,cliff2005} also proposed an open
problem of determining the minimum number of bars among bar tableaux
of a skew shape $\lambda/\mu$. As noted by Clifford
\cite{cliff2003}, this minimum number can be naturally regarded as
the shifted rank, or srank, of $\lambda/\mu$, denoted
$\mathrm{srank}(\lambda/\mu)$. For a skew bar tableau, we present an
algorithm to generate a skew bar tableau without increasing the
number of bars. This algorithm eventually leads to a bar tableau
with the minimum number of bars.

Schur's $Q$-functions arise in the study of the projective
representations of symmetric groups \cite{schur1911}, see also,
Hoffman and Humphreys \cite{hofhum1992}, Humphreys
\cite{humphr1986}, J$\rm{\acute{o}}$zefiak \cite{jozef1989}, Morris
\cite{morri1962, morri1979} and Nazarov \cite{nazar1988}. Shifted
tableaux are closely related to  Schur's $Q$-functions analogous to
the role of ordinary tableaux to the Schur functions. Sagan
\cite{sagan1987} and Worley \cite{worley1984} have independently
developed a combinatorial theory of shifted tableaux, which includes
shifted versions of the Robinson-Schensted-Knuth correspondence,
Knuth's equivalence relations, Sch\"utzenberger's jeu de taquin,
etc. The connections between this combinatorial theory of shifted
tableaux and the theory of projective representations of the
symmetric groups are further explored by Stembridge
\cite{stemb1989}.

Clifford \cite{cliff2005} studied the srank of shifted diagrams for
partitions with distinct parts. Recall that the rank of an ordinary
partition is defined as the number of boxes on the main diagonal of
the corresponding Young diagram. Nazarov and Tarasov
\cite{naztar2002} found an important generalization of the rank of
an ordinary partition to a skew partition in their study of tensor
products of Yangian modules. A general theory of border strip
decompositions and border strip tableaux of skew partitions is
developed by Stanley \cite{stanl2002}, and it has been shown that
the rank of a skew partition is the least number of strips to
construct a minimal border strip decomposition of the skew diagram.
Motivated by Stanley's theorem, Clifford \cite{cliff2005}
generalized the rank of a partition to the rank of a shifted
partition, called srank, in terms of the minimal bar tableaux.

On the other hand, Clifford has noticed that the srank is closely
related to Schur's $Q$-function, as suggested by the work of Stanley
\cite{stanl2002} on the rank of a partition. Stanley introduced a
degree operator by taking the degree of the power sum symmetric
function $p_{\mu}$ as the number of nonzero parts of the indexing
partition $\mu$. Furthermore, Clifford and Stanley \cite{clista2004}
defined the bottom Schur functions to be the sum of the lowest
degree terms in the expansion of the Schur functions in terms of the
power sums. In \cite{cliff2005} Clifford studied the lowest degree
terms in the expansion of Schur's $Q$-functions in terms of power
sum symmetric functions and conjectured that the lowest degree of
the Schur's $Q$-function $Q_{\lambda}$  is equal to the srank of
$\lambda$. Our first result is a proof of this conjecture.

However, in general, the lowest degree of the terms, which appear in
the expansion of the skew Schur's $Q$-function $Q_{\lambda/\mu}$ in
terms of the power sums, is not equal to the srank of the shifted
skew diagram of $\lambda/\mu$. This is different from the case for
ordinary skew partitions and skew Schur functions.  Instead, we will
take an algorithmic approach to the computation of the srank of a
skew partition. It would be interesting to find an algebraic
interpretation in terms of Schur's $Q$-functions.

\section{Shifted diagrams and bar tableaux}\label{sect2}

Throughout this paper we will adopt the notation and terminology on
partitions and symmetric functions in \cite{macdon1995}. A
\emph{partition} $\lambda$ is a weakly decreasing sequence of
positive integers $\lambda_1\geq \lambda_2\geq \ldots\geq
\lambda_k$, denoted $\lambda=(\lambda_1, \lambda_2, \ldots,
\lambda_k)$, and $k$ is called the \emph{length} of $\lambda$,
denoted $\ell(\lambda)$. For convenience we may add sufficient 0's
at the end of $\lambda$ if necessary. If $\sum_{i=1}^k\lambda_i=n$,
we say that $\lambda$ is a partition of the integer $n$, denoted
$\lambda\vdash n$. For each partition $\lambda$ there exists a
geometric representation, known as the Young diagram, which is an
array of squares in the plane justified from the top and left corner
with $\ell(\lambda)$ rows and $\lambda_i$ squares in the $i$-th row.
A partition is said to be \emph{odd} (resp. even) if it has an odd
(resp. even) number of even parts. Let $\mathcal{P}^o(n)$ denote the
set of all partitions of $n$ with only odd parts. We will call a
partition \emph{strict} if all its parts are distinct. Let $\mathcal
{D}(n)$ denote the set of all strict partitions of $n$. For each
partition $\lambda\in \mathcal {D}(n)$, let $S(\lambda)$ be the
shifted diagram of $\lambda$, which is obtained from the Young
diagram by shifting the $i$-th row $(i-1)$ squares to the right for
each $i>1$. For instance, Figure \ref{shifted diagram} illustrates
the shifted diagram of shape $(8,7,5,3,1)$.

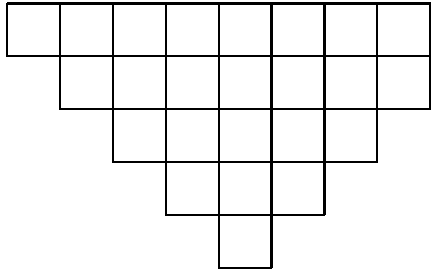
\begin{figure}[h,t]
\setlength{\unitlength}{1pt}
\begin{center}
\begin{picture}(180,120)
\put(10,100){\line(1,0){160}} \put(10,80){\line(1,0){160}}
\put(30,60){\line(1,0){140}}\put(50,40){\line(1,0){100}}
\put(70,20){\line(1,0){60}}\put(90,0){\line(1,0){20}}

\put(10,80){\line(0,1){20}}
\put(30,60){\line(0,1){40}}\put(50,40){\line(0,1){60}}

\put(70,20){\line(0,1){80}}\put(90,0){\line(0,1){100}}
\put(110,0){\line(0,1){100}}\put(130,20){\line(0,1){80}}
\put(150,40){\line(0,1){60}}\put(170,60){\line(0,1){40}}

\end{picture}
\end{center}
\caption{The shifted diagram of shape $(8,7,5,3,1)$}\label{shifted
diagram}
\end{figure}

Given two partitions $\lambda$ and $\mu$, if for each $i$ we have
$\lambda_i\geq \mu_i$, then the skew partition $\lambda/\mu$ is
defined to be the diagram obtained from the diagram of $\lambda$
by removing the diagram of $\mu$ at the top-left corner.
Similarly, the skew shifted diagram $S(\lambda/\mu)$ is defined as
the set-theoretic difference of $S(\lambda)$ and $S(\mu)$.

Now we recall the definitions of bars and bar tableaux as given in
Hoffman and Humphreys \cite{hofhum1992}. Let $\lambda\in \mathcal
{D}(n)$ be a partition with length $\ell(\lambda)=k$. Fixing an odd
positive integer $r$, three subsets $I_{+}, I_{0}, I_{-}$ of
integers between $1$ and $k$ are defined as follows:
\begin{eqnarray*}
I_{+}& = &\{i: \lambda_{j+1}<\lambda_i-r<\lambda_j\: \mbox{for
some }
j\leq k,\: \mbox {taking}\:\lambda_{k+1}=0\},\\[5pt]
I_{0} & = & \{i: \lambda_i=r\},\\[5pt]
I_{-} & = & \{i: r-\lambda_{i}=\lambda_j \:\mbox{for some} \:j\:
\mbox{with} \:i<j\leq k\}.
\end{eqnarray*}

Let $I(\lambda,r)=I_{+}\cup I_{0}\cup I_{-}$. For each $i\in
I(\lambda,r)$, we define a new strict partition $\lambda(i,r)$ of
$\mathcal {D}(n-r)$ in the following way:
\begin{itemize}
\item[(1)] If $i\in I_{+}$, then $\lambda_i>r$, and let $\lambda(i,r)$ be
the partition obtained from $\lambda$ by removing $\lambda_i$ and
inserting $\lambda_i-r$ between $\lambda_j$ and $\lambda_{j+1}$.
\item[(2)] If $i\in I_{0}$, let $\lambda(i,r)$ be the partition obtained
from $\lambda$ by removing $\lambda_i$.
\item[(3)] If $i\in I_{-}$, then let $\lambda(i,r)$ be the partition obtained
from $\lambda$ by removing both $\lambda_i$ and $\lambda_j$.
\end{itemize}

Meanwhile, for each $i\in I(\lambda,r)$, the associated $r$-bar is
given as follows:
\begin{itemize}
\item[(1')] If $i\in I_{+}$, the $r$-bar consists of the rightmost $r$
squares in the $i$-th row of $S(\lambda)$, and we say that the
$r$-bar is of Type $1$.
\item[(2')] If $i\in I_{0}$, the $r$-bar
consists of all the squares of the $i$-th row of $S(\lambda)$, and
we say that the $r$-bar is of Type $2$.
\item[(3')] If $i\in I_{-}$, the $r$-bar consists of all the squares of the
$i$-th and $j$-th rows, and we say that the $r$-bar is of Type $3$.
\end{itemize}
For example, as shown in Figure \ref{bar tableau}, the squares
filled with $6$ are a $7$-bar of Type $1$, the squares filled with
$4$ are a $3$-bar of Type $2$, and the squares filled with $3$ are a
$7$-bar of Type $3$.

\begin{figure}[h,t]
\setlength{\unitlength}{1pt}
\begin{center}
\begin{picture}(180,103)
\put(10,100){\line(1,0){180}} \put(10,80){\line(1,0){180}}
\put(30,60){\line(1,0){140}}\put(50,40){\line(1,0){120}}
\put(70,20){\line(1,0){60}}\put(90,0){\line(1,0){20}}

\put(10,80){\line(0,1){20}}
\put(30,60){\line(0,1){40}}\put(50,40){\line(0,1){60}}

\put(70,20){\line(0,1){80}}\put(90,0){\line(0,1){100}}
\put(110,0){\line(0,1){100}}\put(130,20){\line(0,1){80}}
\put(150,40){\line(0,1){60}}\put(170,40){\line(0,1){60}}
\put(190,80){\line(0,1){20}}

\put(18,87){$1$}\put(38,87){$1$}\put(58,87){$6$}
\put(78,87){$6$}\put(98,87){$6$}\put(118,87){$6$}
\put(138,87){$6$}\put(158,87){$6$}\put(178,87){$6$}

\put(38,67){$1$}\put(58,67){$2$}\put(78,67){$2$}
\put(98,67){$2$}\put(118,67){$5$}\put(138,67){$5$}
\put(158,67){$5$}

\put(58,47){$3$}\put(78,47){$3$}\put(98,47){$3$}
\put(118,47){$3$}\put(138,47){$3$}\put(158,47){$3$}

\put(78,27){$4$}\put(98,27){$4$}\put(118,27){$4$}

\put(98,7){$3$}
\end{picture}
\end{center}
\caption{A bar tableau of shape $(9,7,6,3,1)$}\label{bar tableau}
\end{figure}

A \emph{bar tableau} of shape $\lambda$ is an array of positive
integers of shape $S(\lambda)$ subject to the following conditions:
\begin{itemize}
\item[(1)] It is weakly increasing in every row;

\item[(2)] The number of parts equal to $i$ is odd for each positive integer
$i$;

\item[(3)] Each positive integer $i$ can appear in at most two rows,
and if $i$ appears in two rows, then these two rows must begin with
$i$;

\item[(4)] The composition obtained by removing all squares
filled with integers larger than some $i$ has distinct parts.
\end{itemize}

We say that a bar tableau $T$ is of type
$\rho=(\rho_1,\rho_2,\ldots)$ if the total number of $i$'s appearing
in $T$ is  $\rho_i$. For example, the bar tableau in Figure
\ref{weight} is of type $(3,1,1,1)$. For a bar tableau $T$ of shape
$\lambda$,  we define its weight $wt(T)$ recursively by the
following procedure. If $T$ is empty, let $wt(T)=1$. Let
$\varepsilon(\lambda)$ denote the parity of the partition $\lambda$,
i.e., $\varepsilon(\lambda)=0$ if $\lambda$ has an even number of
even parts; otherwise, $\varepsilon(\lambda)=1$. Suppose that the
largest numbers in $T$ form an $r$-bar, which is associated with an
index $i\in I(\lambda, r)$. Let $j$ be the integer that occurrs in
the definitions of $I_{+}$ and $I_{-}$. Let $T'$ be the bar tableau
of shape $\lambda(i, r)$ obtained from $T$ by removing this $r$-bar.
Now, let
\begin{equation}
wt(T)=n_i\, wt(T'),
\end{equation}
where
\begin{equation}
n_i=\left\{\begin{array}{cc}
(-1)^{j-i}2^{1-\varepsilon(\lambda)},&
\mbox{if}\  i\in I_{+},\\[6pt]
(-1)^{\ell(\lambda)-i},& \mbox{if}\  i\in I_{0},\\[6pt]
(-1)^{j-i+\lambda_i}2^{1-\varepsilon(\lambda)},& \mbox{if}\  i\in
I_{-}.
\end{array}
\right.
\end{equation}
For example, the weight of the bar tableau $T$ in Figure
\ref{weight} equals
\begin{equation}
wt(T)=(-1)^{1-1}2^{1-0}\cdot(-1)^{1-1}2^{1-1}\cdot(-1)^{2-2}
\cdot(-1)^{1-1}=2.
\end{equation}

\begin{figure}[h,t]
\setlength{\unitlength}{1pt}
\begin{center}
\begin{picture}(180,40)
\put(40,40){\line(1,0){100}}\put(40,20){\line(1,0){100}}
\put(60,0){\line(1,0){20}}

\put(40,20){\line(0,1){20}}\put(60,0){\line(0,1){40}}
\put(80,0){\line(0,1){40}}\put(100,20){\line(0,1){20}}
\put(120,20){\line(0,1){20}}\put(140,20){\line(0,1){20}}
\put(48,26){$1$}\put(68,26){$1$}
\put(88,26){$1$}\put(108,26){$3$}\put(128,26){$4$} \put(68,6){$2$}
\end{picture}
\end{center}
\caption{A bar tableau of type $(3,1,1,1)$}\label{weight}
\end{figure}

The following lemma will be used  in Section 3 to determine whether
certain terms will vanish in the power sum expansion of Schur's
$Q$-functions indexed by partitions with two distinct parts.

\begin{lem} \label{vanishbar}
Let $\lambda=(\lambda_1,\lambda_2)$ be a strict partition with the
two parts $\lambda_1$ and $\lambda_2$ having the same parity. Given
an partition $\sigma=(\sigma_1,\sigma_2)\in
\mathcal{P}^o(|\lambda|)$, if $\sigma_2<\lambda_2$, then among all
bar tableaux of shape $\lambda$ there exist only two bar tableaux of
type $\sigma$, say $T_1$ and $T_2$, and furthermore, we have
$wt(T_1)+wt(T_2)=0$.
\end{lem}

\pf Suppose that both $\lambda_1$ and $\lambda_2$ are even. The case
when $\lambda_1$ and $\lambda_2$ are odd numbers can be proved
similarly. Note that $\sigma_2<\lambda_{2}<\lambda_{1}$. By putting
$2$'s in the last $\sigma_2$ squares of the second row and then
filling the remaining squares in the diagram with $1$'s, we obtain
one tableau $T_1$. By putting $2$'s in the last $\sigma_2$ squares
of the first row and then filling the remaining squares with $1$'s,
we obtain another tableau $T_2$. Clearly, both $T_1$ and $T_2$ are
bar tableaux of shape $\lambda$ and type $\sigma$, and they are the
only two such bar tableaux. We notice that
\begin{equation}
wt(T_1)=(-1)^{2-2}2^{1-0}\cdot (-1)^{2-1+\lambda_1} 2^{1-1}=-2.
\end{equation}
While, for the weight of $T_2$, there are two cases to consider.
If $\lambda_1-\sigma_2>\lambda_2$, then
\begin{equation}
wt(T_2)=(-1)^{1-1}2^{1-0}\cdot
(-1)^{2-1+\lambda_1-\sigma_2}2^{1-1}=2.
\end{equation}
If $\lambda_1-\sigma_2<\lambda_2$, then
\begin{equation}
wt(T_2)=(-1)^{2-1}2^{1-0}\cdot (-1)^{2-1+\lambda_2}2^{1-1}=2.
\end{equation}
Thus we have $wt(T_2)=2$ in either case, so the relation
$wt(T_1)+wt(T_2)=0$ holds. \qed

For example, taking $\lambda=(8,6)$ and $\sigma=(11,3)$, the two bar
tableaux $T_1$ and $T_2$ in the above lemma are depicted  as in
Figure \ref{2-bar tableaux1}.

\begin{figure}[h,t]
\setlength{\unitlength}{1.8pt}
\begin{center}
\begin{picture}(180,40)
\put(0,30){\line(1,0){80}}\put(0,20){\line(1,0){80}}
\put(10,10){\line(1,0){60}}\put(0,30){\line(0,-1){10}}
\multiput(10,30)(10,0){7}{\line(0,-1){20}}
\put(80,30){\line(0,-1){10}} \multiput(4,23)(10,0){8}{$1$}
\multiput(14,13)(10,0){3}{$1$} \multiput(44,13)(10,0){3}{$2$}
\put(40,0){$T_1$}
\put(110,30){\line(1,0){80}}\put(110,20){\line(1,0){80}}
\put(120,10){\line(1,0){60}}\put(110,30){\line(0,-1){10}}
\multiput(120,30)(10,0){7}{\line(0,-1){20}}
\put(190,30){\line(0,-1){10}} \multiput(114,23)(10,0){5}{$1$}
\multiput(164,23)(10,0){3}{$2$}\multiput(124,13)(10,0){6}{$1$}
\put(150,0){$T_2$}
\end{picture}
\end{center}
\caption{Two bar tableaux of shape $(8,6)$ and type
$(11,3)$}\label{2-bar tableaux1}
\end{figure}

Clifford gave a natural generalization of bar tableaux to skew
shapes \cite{cliff2005}. Formally, a \emph{skew bar tableau} of
shape $\lambda/\mu$ is an assignment of nonnegative integers to the
squares of $S(\lambda)$ such that in addition to the above four
conditions (1)-(4) we further impose the condition that
\begin{itemize}
\item [(5)] the partition obtained by removing all squares
filled with positive integers and reordering the remaining rows is
$\mu$.
\end{itemize}

For example, taking the skew partition $(8,6,5,4,1)/(8,2,1)$, Figure
\ref{skew bar tableau} is a skew bar tableau of such shape.

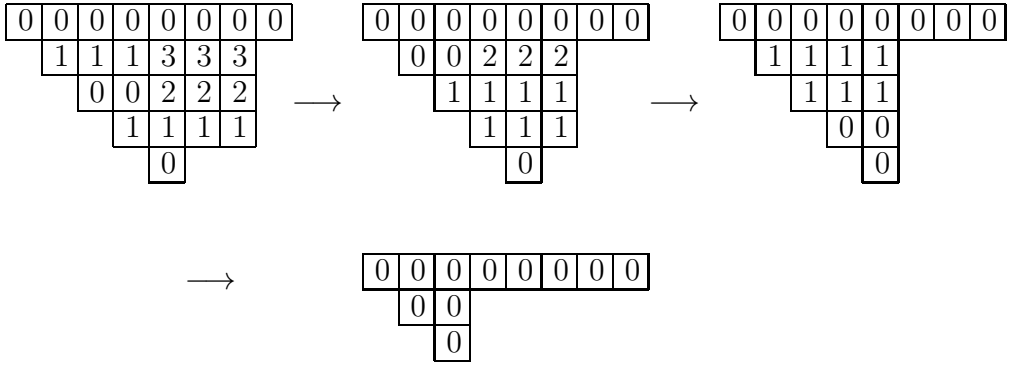
\begin{figure}[h,t]
\setlength{\unitlength}{0.9pt}
\begin{center}
\begin{picture}(180,160)
\put(-125,150){\line(1,0){120}}\put(-125,135){\line(1,0){120}}
\put(-110,120){\line(1,0){90}}\put(-95,105){\line(1,0){75}}
\put(-80,90){\line(1,0){60}}\put(-65,75){\line(1,0){15}}
\put(-125,135){\line(0,1){15}}\put(-110,120){\line(0,1){30}}
\put(-95,105){\line(0,1){45}}\put(-80,90){\line(0,1){60}}
\put(-65,75){\line(0,1){75}}\put(-50,75){\line(0,1){75}}
\put(-35,90){\line(0,1){60}}\put(-20,90){\line(0,1){60}}
\put(-5,135){\line(0,1){15}} \put(-5,105){$\longrightarrow$}

\put(-120,139){$0$}\put(-105,139){$0$}\put(-90,139){$0$}
\put(-75,139){$0$}\put(-60,139){$0$}\put(-45,139){$0$}
\put(-30,139){$0$}\put(-15,139){$0$}

\put(-105,124){$1$}\put(-90,124){$1$}
\put(-75,124){$1$}\put(-60,124){$3$}\put(-45,124){$3$}
\put(-30,124){$3$}

\put(-90,109){$0$}
\put(-75,109){$0$}\put(-60,109){$2$}\put(-45,109){$2$}
\put(-30,109){$2$}

\put(-75,94){$1$}\put(-60,94){$1$}\put(-45,94){$1$}
\put(-30,94){$1$}
 \put(-60,79){$0$}

\put(25,150){\line(1,0){120}}\put(25,135){\line(1,0){120}}
\put(40,120){\line(1,0){75}}\put(55,105){\line(1,0){60}}
\put(70,90){\line(1,0){45}}\put(85,75){\line(1,0){15}}
\put(25,135){\line(0,1){15}}\put(40,120){\line(0,1){30}}
\put(55,105){\line(0,1){45}}\put(70,90){\line(0,1){60}}
\put(85,75){\line(0,1){75}}\put(100,75){\line(0,1){75}}
\put(115,90){\line(0,1){60}}\put(130,135){\line(0,1){15}}
\put(145,135){\line(0,1){15}}\put(145,105){$\longrightarrow$}

\put(30,139){$0$}\put(45,139){$0$}\put(60,139){$0$}
\put(75,139){$0$}\put(90,139){$0$}\put(105,139){$0$}
\put(120,139){$0$}\put(135,139){$0$}

\put(45,124){$0$}\put(60,124){$0$}
\put(75,124){$2$}\put(90,124){$2$}\put(105,124){$2$}

\put(60,109){$1$}
\put(75,109){$1$}\put(90,109){$1$}\put(105,109){$1$}

\put(75,94){$1$}\put(90,94){$1$}\put(105,94){$1$}

\put(90,79){$0$}

\put(175,150){\line(1,0){120}}\put(175,135){\line(1,0){120}}
\put(190,120){\line(1,0){60}}\put(205,105){\line(1,0){45}}
\put(220,90){\line(1,0){30}}\put(235,75){\line(1,0){15}}

\put(175,135){\line(0,1){15}}\put(190,120){\line(0,1){30}}
\put(205,105){\line(0,1){45}}\put(220,90){\line(0,1){60}}
\put(235,75){\line(0,1){75}}\put(250,75){\line(0,1){75}}
\put(265,135){\line(0,1){15}}\put(280,135){\line(0,1){15}}
\put(295,135){\line(0,1){15}}

\put(180,139){$0$}\put(195,139){$0$}\put(210,139){$0$}
\put(225,139){$0$}\put(240,139){$0$}\put(255,139){$0$}
\put(270,139){$0$}\put(285,139){$0$}

\put(195,124){$1$}\put(210,124){$1$}
\put(225,124){$1$}\put(240,124){$1$}

\put(210,109){$1$} \put(225,109){$1$}\put(240,109){$1$}

\put(225,94){$0$}\put(240,94){$0$}

\put(240,79){$0$}

\put(-50,30){$\longrightarrow$}
\put(25,45){\line(1,0){120}}\put(25,30){\line(1,0){120}}
\put(40,15){\line(1,0){30}}\put(55,0){\line(1,0){15}}

\put(25,30){\line(0,1){15}}\put(40,15){\line(0,1){30}}
\put(55,0){\line(0,1){45}}\put(70,0){\line(0,1){45}}
\put(85,30){\line(0,1){15}}\put(100,30){\line(0,1){15}}
\put(115,30){\line(0,1){15}}\put(130,30){\line(0,1){15}}
\put(145,30){\line(0,1){15}}

\put(30,34){$0$}\put(45,34){$0$}\put(60,34){$0$}
\put(75,34){$0$}\put(90,34){$0$}\put(105,34){$0$}
\put(120,34){$0$}\put(135,34){$0$}

 \put(45,19){$0$}\put(60,19){$0$}
\put(60,4){$0$}

\end{picture}
\end{center}
\caption{Checking the legality of a skew bar tableau} \label{skew
bar tableau}
\end{figure}

A bar tableau of shape $\lambda$ is said to be \emph{minimal} if
there does not exist a bar tableau with fewer bars. Motivated by
Stanley's results in \cite{stanl2002}, Clifford defined the srank of
a shifted partition $S(\lambda)$, denoted ${\rm srank}(\lambda)$, as
the number of bars in a minimal bar tableau of shape $\lambda$
\cite{cliff2005}. Clifford also gave the following formula for ${\rm
srank}(\lambda)$.

\begin{thm}[{\cite[Theorem 4.1]{cliff2005}}]\label{min bar}
Given a strict partition $\lambda$, let $o$ be the number of odd
parts of $\lambda$, and let $e$ be the number of even parts. Then
${\rm srank}(\lambda)=\max(o,e+(\ell(\lambda) \ \mathrm{mod}\  2))$.
\end{thm}

Next we consider the number of bars in a minimal skew bar tableau of
shape $\lambda/\mu$. Note that the squares filled with $0$'s in the
skew bar tableau give rise to a shifted diagram of shape $\mu$ by
reordering the rows. Let $o_ r$ (resp. $e_r$) be the number of
nonempty rows of odd (resp. even) length with blank squares, and let
$o_s$ (resp. $e_s$) be the number of rows of $\lambda$ with some
squares filled with $0$'s and  an odd (resp. even) number of blank
squares. It is obvious that the number of bars in a minimal skew bar
tableau is greater than or equal to
$$o_s+2e_s+\max(o_r,e_r+((e_r+o_r)\ \mathrm{mod}\ 2)).$$
In fact the above quantity has been considered by Clifford
{\cite{cliff2003}}. Observe  that this quantity depends on the
positions of the 0's.

It should be remarked that a legal bar tableau of shape
$\lambda/\mu$ may not exist once the positions of $0$'s are fixed.
One open problem proposed by Clifford \cite{cliff2003} is to find a
characterization of ${\rm srank}(\lambda/\mu)$. In Section 5 we will
give an algorithm to compute the srank of a skew shape.

\section{Clifford's conjecture}\label{sect3}

In this section, we aim to show  that the lowest degree of the power
sum expansion of a Schur's $Q$-function $Q_{\lambda}$ equals ${\rm
srank}(\lambda)$. Let us recall relevant terminology on Schur's
$Q$-functions.  Let $x=(x_1,x_2,\ldots)$ be an infinite sequence of
independent indeterminates. We define the symmetric functions
$q_k=q_k(x)$ in $x_1,x_2,\ldots$ for all integers $k$ by the
following expansion of the formal power series in $t$:
$$\prod_{i\geq 1}\frac{1+x_it}{1-x_it}=\sum_{k}q_{k}(x)t^k.$$
In particular, $q_k=0$ for $k<0$ and $q_0=1$. It immediately follows
that
\begin{equation}\label{eq-def}
\sum_{i+j=n}(-1)^iq_iq_j=0,
\end{equation}
for all $n\geq 1$.
Let $Q_{(a)}=q_a$ and
$$Q_{(a,b)}=q_aq_b+2\sum_{m=1}^b(-1)^m q_{a+m}q_{b-m}.$$
From \eqref{eq-def} we see that $Q_{(a,b)}=-Q_{(b,a)}$ and thus
$Q_{(a,a)}=0$ for any $a,b$. In general, for any strict partition
$\lambda$, the symmetric function $Q_{\lambda}$ is defined by the
recurrence relations:
\begin{eqnarray}
Q_{(\lambda_1,\ldots,\lambda_{2k+1})}&=& \sum_{m=1}^{2k+1}
(-1)^{m+1}
q_{\lambda_m}Q_{(\lambda_1,\ldots,\hat{\lambda}_m,\ldots,\lambda_{2k+1})},\\[5pt]
Q_{(\lambda_1,\ldots,\lambda_{2k})}&=& \sum_{m=2}^{2k} (-1)^{m}
Q_{(\lambda_1,\lambda_m)}Q_{(\lambda_2,\ldots,\hat{\lambda}_m,\ldots,\lambda_{2k})},
\end{eqnarray}
where $\hat{}$ stands for a missing entry.

It was known that $Q_{\lambda}$ can be  also defined as the
specialization at $t=-1$ of the Hall-Littlewood functions associated
with $\lambda$ \cite{macdon1995}. Originally, these $Q_{\lambda}$
symmetric functions were introduced in order to express irreducible
projective characters of the symmetric groups \cite{schur1911}. Note
that the irreducible projective representations of $S_n$ are in
one-to-one correspondence with partitions of $n$ with distinct
parts, see \cite{jozef1989, stemb1988, stemb1989}. For any
$\lambda\in \mathcal{D}(n)$, let $\langle\lambda\rangle$ denote the
character of the irreducible projective or spin representation
indexed by $\lambda$.  Morris \cite{morri1965} has found a
combinatorial rule for calculating the characters, which is the
projective analogue of the Murnaghan-Nakayama rule. In terms of bar
tableaux, Morris's theorem reads as follows:

\begin{thm}[\cite{morri1965}]\label{mnrule}
Let $\lambda\in \mathcal{D}(n)$ and $\pi\in
\mathcal{P}^o(n)$. Then
\begin{equation}\label{mnruleeq}
\langle\lambda\rangle(\pi)=\sum_{T}wt(T)
\end{equation}
where the sum ranges over all bar tableaux of shape $\lambda$ and
type $\pi$.
\end{thm}

The above theorem for projective characters implies the following
formula, which will be used later in the proof of Lemma \ref{len2}.

\begin{cor} \label{2odd} Let $\lambda$ be a strict partition of length $2$.
Suppose that the two parts $\lambda_1,\lambda_2$ are both odd. Then
we have
\begin{equation}
\langle\lambda\rangle(\lambda)=-1.
\end{equation}
\end{cor}

\pf Let $T$ be the bar tableau obtained by filling the last
$\lambda_2$ squares in the first row of $S(\lambda)$ with $2$'s and
the remaining squares with $1$'s, and let $T'$ be the bar tableau
obtained by filling the first row of $S(\lambda)$ with $1$'s and the
second row with $2$'s. Clearly, $T$ and $T'$ are of the same type
$\lambda$. Let us first consider the weight of $T$. If
$\lambda_1-\lambda_2<\lambda_2$, then
$$
wt(T)=(-1)^{2-1} 2^{1-0}\cdot (-1)^{2-1+\lambda_2}2^{1-1}=-2.$$ If
$\lambda_1-\lambda_2>\lambda_2$, then
$$
 wt(T)=(-1)^{1-1} 2^{1-0}\cdot
(-1)^{2-1+\lambda_1-\lambda_2}2^{1-1}=-2.
$$
In both cases, the weight of $T'$ equals
$$wt(T')=(-1)^{2-2}\cdot (-1)^{1-1}=1.$$
Since there are only two bar tableaux, $T$ and $T'$, of type
$\lambda$, the  corollary immediately follows  from Theorem
\ref{mnrule}. \qed

Let $p_k(x)$ denote the $k$-th power sum symmetric functions, i.e.,
$p_k(x)=\sum_{i\geq 1}x_i^k$. For any partition
$\lambda=(\lambda_1,\lambda_2,\cdots)$, let
$p_{\lambda}=p_{\lambda_1}p_{\lambda_2}\cdots$. The fundamental
connection between $Q_{\lambda}$ symmetric functions and the
projective representations of the symmetric group is as follows.

\begin{thm}[\cite{schur1911}]\label{conn}
Let $\lambda\in \mathcal{D}(n)$. Then we have
\begin{equation}
Q_{\lambda}=\sum_{\pi\in \mathcal{P}^o(n)}
2^{[\ell(\lambda)+\ell(\pi)+\varepsilon(\lambda)]/2}
\langle\lambda\rangle(\pi)\frac{p_{\pi}}{z_{\pi}},
\end{equation}
where
$$z_{\pi}=1^{m_1}m_1!\cdot 2^{m_2}m_2!\cdot \cdots, \quad \mbox{if $\pi=\langle 1^{m_1}2^{m_2}\cdots \rangle$.}$$
\end{thm}

Stanley \cite{stanl2002} introduced a degree operator on symmetric
functions by defining $\deg(p_i)=1$, and so
$\deg(p_{\nu})=\ell(\nu)$. Clifford \cite{cliff2005} applied this
operator to Schur's $Q$-functions and obtained the following lower
bound from Theorem \ref{conn}.
\begin{cor}[{\cite[Corollary 6.2]{cliff2005}}] \label{atleast}
The terms of the lowest degree in $Q_{\lambda}$ have degree at least
${\rm srank}(\lambda)$.
\end{cor}

 The following conjecture is proposed by Clifford:

\begin{con}[{\cite[Conjecture 6.4]{cliff2005}}]
The terms of the lowest degree in $Q_{\lambda}$ have degree  ${\rm
srank}(\lambda)$.
\end{con}

Our proof of the above conjecture depends on the Pfaffian formula
for Schur's $Q$-functions. Given a skew-symmetric matrix
$A=(a_{i,j})$ of even size $2n\times 2n$, the \emph{Pfaffian} of
$A$, denoted {\rm Pf}(A), is defined by
$${\rm Pf}(A)=\sum_{\pi}(-1)^{{\rm cr}(\pi)} a_{i_1j_1}\cdots a_{i_nj_n},$$
where the sum ranges over all set partitions $\pi$ of $\{1,2,\cdots,
2n\}$ into two element blocks $i_k<j_k$ and $cr(\pi)$ is the number
of crossings of $\pi$, i.e., the number of pairs $h<k$ for which
$i_h<i_k<j_h<j_k$.

\begin{thm}[\cite{macdon1995}] \label{pfexp}
Given a strict partition
$\lambda=(\lambda_1,\lambda_2,\ldots,\lambda_{2n})$ satisfying
$\lambda_1>\ldots>\lambda_{2n}\geq 0$, let
$M_{\lambda}=(Q_{(\lambda_i,\lambda_j)})$.
Then we have
$$Q_{\lambda}={\rm Pf}(M_{\lambda}).$$
\end{thm}

We first prove that Clifford's conjecture holds for strict
partitions of length less than three. The proof for the general case
relies on this special case.

\begin{lem}\label{len2}
Let $\lambda$ be a strict partition of length $\ell(\lambda)<3$.
Then the terms of the lowest degree in $Q_{\lambda}$ have degree
${\rm srank}(\lambda)$.
\end{lem}
\pf In view of Theorem \ref{mnrule} and Theorem \ref{conn}, if there
exists a unique bar tableau of shape $\lambda$ and type $\pi$, then
the coefficient of $p_{\pi}$ is nonzero in the expansion of
$Q_{\lambda}$. There are five cases to consider.
\begin{itemize}

\item[(1)] $\ell(\lambda)=1$ and $\lambda_1$ is
odd. Clearly, we have ${\rm srank}(\lambda)=1$. Note that there
exists a unique bar tableau $T$ of shape $\lambda$ and of type
$\lambda$ with all squares of $S(\lambda)$ filled with $1$'s.
Therefore, the coefficient of $p_{\lambda}$ in the power sum
expansion of $Q_{\lambda}$ is nonzero and the lowest degree of
$Q_{\lambda}$ is $1$.

\item[(2)] $\ell(\lambda)=1$ and $\lambda_1$ is
even. We see that ${\rm srank}(\lambda)=2$. Since the bars are all
of odd size, there does not exist any bar tableau of shape $\lambda$
and of type $\lambda$. But there is a unique bar tableau $T$ of
shape $\lambda$ and of type $(\lambda_1-1,1)$, which is obtained by
filling the rightmost square of $S(\lambda)$ with $2$ and the
remaining squares with $1$'s. So the coefficient of
$p_{(\lambda_1-1,1)}$ in the power sum expansion of $Q_{\lambda}$ is
nonzero and the terms of the lowest degree in $Q_{\lambda}$ have
degree $2$.

\item[(3)] $\ell(\lambda)=2$ and the two parts
$\lambda_1,\lambda_2$ have different parity. In this case, we have
${\rm srank}(\lambda)=1$. Note that there exists a unique bar
tableau $T$ of shape $\lambda$ and of type $(\lambda_1+\lambda_2)$,
which is obtained by filling all the squares of $S(\lambda)$ with
$1$'s. Thus, the coefficient of $p_{\lambda_1+\lambda_2}$ in the
power sum expansion of $Q_{\lambda}$ is nonzero and the terms of
lowest degree in $Q_{\lambda}$ have degree $1$.

\item[(4)] $\ell(\lambda)=2$ and the two parts
$\lambda_1,\lambda_2$ are both even. It is easy to see that ${\rm
srank}(\lambda)=2$. Since there exists a unique bar tableau $T$ of
shape $\lambda$ and of type $(\lambda_1-1,\lambda_2+1)$, which is
obtained by filling the rightmost $\lambda_2+1$ squares in the first
row of $S(\lambda)$ with $2$'s and the remaining squares with $1$'s,
the coefficient of $p_{(\lambda_1-1,\lambda_2+1)}$ in the power sum
expansion of $Q_{\lambda}$ is nonzero; hence the lowest degree of
$Q_{\lambda}$ is equal to $2$.

\item[(5)] $\ell(\lambda)=2$ and the two parts
$\lambda_1,\lambda_2$ are both odd. In this case,  we have ${\rm
srank}(\lambda)=2$. By Corollary \ref{2odd},  the coefficient of
$p_{\lambda}$ in the power sum expansion of $Q_{\lambda}$ is
nonzero, and therefore the terms of the lowest degree in
$Q_{\lambda}$ have degree $2$.
\end{itemize}
This completes the proof. \qed

Given a strict partition $\lambda$, we consider the Pfaffian
expansion of $Q_{\lambda}$ as shown in Theorem \ref{pfexp}. To prove
Clifford's conjecture, we need to determine which terms may appear
in the expansion of $Q_{\lambda}$ in terms of power sum symmetric
functions. Suppose that the Pfaffian expansion of $Q_{\lambda}$ is
as follows:
\begin{equation}\label{q-expand}
{\rm Pf}(M_{\lambda})=\sum_{\pi}(-1)^{{\rm cr}(\pi)}
Q_{(\lambda_{\pi_1},\lambda_{\pi_2})}\cdots
Q_{(\lambda_{\pi_{2m-1}},\lambda_{\pi_{2m}})},
\end{equation}
where the sum ranges over all set partitions $\pi$ of $\{1,2,\cdots,
2m\}$ into two element blocks
$\{(\pi_1,\pi_2),\ldots,(\pi_{2m-1},\pi_{2m})\}$ with
$\pi_1<\pi_3<\cdots<\pi_{2m-1}$ and $\pi_{2k-1}<\pi_{2k}$ for any
$k$. For the above expansion of $Q_{\lambda}$, the following two
lemmas will be used to choose certain lowest degree terms in the
power sum expansion of $Q_{(\lambda_i,\lambda_j)}$ in the matrix
$M_\lambda$.

\begin{lem}\label{lemma1} Suppose that $\lambda$ has both odd parts and even
parts. Let $\lambda_{i_1}$ (resp. $\lambda_{j_1}$) be the largest
odd (resp. even) part of $\lambda$. If the power sum symmetric
function $p_{\lambda_{i_1}+\lambda_{j_1}}$ appears in the terms of
lowest degree originated from the product
$Q_{(\lambda_{\pi_1},\lambda_{\pi_2})}\cdots
Q_{(\lambda_{\pi_{2m-1}},\lambda_{\pi_{2m}})}$ as in the expansion
\eqref{q-expand}, then we have $(\pi_1,\pi_2)=(i_1,j_1)$.
\end{lem}

\pf Without loss of generality, we may assume that $\lambda_{i_1}>
\lambda_{j_1}$. By Lemma \ref{len2}, the term
$p_{\lambda_{i_1}+\lambda_{j_1}}$ appears in $Q_{(\lambda_{i_1},
\lambda_{j_1})}$ with nonzero coefficients. Since $\lambda_{i_1},
\lambda_{j_1}$ are the largest odd and even parts,
$p_{\lambda_{i_1}+\lambda_{j_1}}$ does not appear as a factor of any
term of the lowest degree in the expansion of $Q_{(\lambda_{i_k},
\lambda_{j_k})}$, where $\lambda_{i_k}$ and $\lambda_{j_k}$ have
different parity. Meanwhile, if $\lambda_{i_k}$ and $\lambda_{j_k}$
have the same parity, then we consider the bar tableaux of shape
$(\lambda_{i_k}, \lambda_{j_k})$ and of type
$(\lambda_{i_1}+\lambda_{j_1}, \lambda_{i_k}+
\lambda_{j_k}-\lambda_{i_1}-\lambda_{j_1})$. Observe that
$\lambda_{i_k}+
\lambda_{j_k}-\lambda_{i_1}-\lambda_{j_1}<\lambda_{j_k}$. Since the
lowest degree of $Q_{(\lambda_{i_k}, \lambda_{j_k})}$ is $2$, from
Lemma \ref{vanishbar} it follows that
$p_{\lambda_{i_1}+\lambda_{j_1}}$ can not be a factor of any term of
lowest degree in the power sum expansion of $Q_{(\lambda_{i_k},
\lambda_{j_k})}$. This completes the proof. \qed

\begin{lem}\label{lemma2}
Suppose that $\lambda$ only has even parts. Let $\lambda_1,
\lambda_2$ be the two largest parts of $\lambda$ (allowing
$\lambda_2=0$). If the power sums $p_{\lambda_1-1}p_{\lambda_2+1}$
appears in the terms of the lowest degree given by the product
$Q_{(\lambda_{\pi_1},\lambda_{\pi_2})}\cdots
Q_{(\lambda_{\pi_{2m-1}},\lambda_{\pi_{2m}})}$ as in
\eqref{q-expand}, then we have $(\pi_1,\pi_2)=(1,2)$.
\end{lem}

\pf From Case (4) of the proof of Lemma \ref{len2} it follows that
$p_{\lambda_1-1}p_{\lambda_2+1}$ appears as a term of the lowest
degree in the power sum expansion of $Q_{(\lambda_1,\lambda_2)}$. We
next consider the power sum expansion of any other
$Q_{(\lambda_i,\lambda_j)}$. First, we consider the case when
$\lambda_i+\lambda_j>\lambda_2+1$ and $\lambda_i \leq\lambda_2$.
Since $\lambda_i+\lambda_j-(\lambda_2+1)<\lambda_j$, by Lemma
\ref{vanishbar}, the term $p_{\lambda_2+1}$  is not a factor of any
term of the lowest degree in the power sum expansion of
$Q_{(\lambda_i,\lambda_j)}$. Now we are left with the case when
$\lambda_i+\lambda_j>\lambda_1-1$ and $\lambda_i\leq \lambda_1-2$.
Since $\lambda_i+\lambda_j-(\lambda_1-1)<\lambda_j$, by Lemma
\ref{vanishbar} the term $p_{\lambda_1-1}$ does not appear as a
factor in the terms of the lowest degree of
$Q_{(\lambda_i,\lambda_j)}$.  So we have shown that if either
$p_{\lambda_2+1}$ or $p_{\lambda_1-1}$ appears as a factor of some
lowest degree term for $Q_{(\lambda_i,\lambda_j)}$, then we deduce
that $\lambda_i=\lambda_1$. Moreover, if  both $p_{\lambda_1-1}$ and
$p_{\lambda_2+1}$ are factors of the lowest degree terms in the
power sum expansion of $Q_{(\lambda_1,\lambda_j)}$, then we have
$\lambda_j=\lambda_2$. The proof is complete. \qed

We now present the main result of this paper.

\begin{thm}
For any $\lambda\in\mathcal{D}(n)$, the terms of the lowest degree
in $Q_\lambda$ have degree ${\rm srank}(\lambda)$.
\end{thm}

\pf We  write the strict partition $\lambda$ in the form
$(\lambda_1,\lambda_2,\ldots,\lambda_{2m})$, where
$\lambda_1>\ldots>\lambda_{2m}\geq 0$. Suppose that the partition
$\lambda$ has $o$ odd parts and $e$ even parts (including $0$ as a
part). For the sake of presentation, let
$(\lambda_{i_1},\lambda_{i_2},\ldots,\lambda_{i_o})$ denote the
sequence of odd parts in decreasing order, and let
$(\lambda_{j_1},\lambda_{j_2},\ldots,\lambda_{j_e})$ denote the
sequence of even parts in decreasing order.

We first consider the case  $o\geq e$. In this case, it will be
shown that  ${\rm srank}(\lambda)=o$. By Theorem \ref{min bar}, if
$\lambda_{2m}>0$, i.e., $\ell(\lambda)=2m$, then we have
$${\rm srank}(\lambda)=\max(o,e+0)=o.$$ If $\lambda_{2m}=0$, i.e.,
$\ell(\lambda)=2m-1$, then we still have
$${\rm srank}(\lambda)=\max(o,(e-1)+1)=o.$$

Let
$$A=p_{\lambda_{i_1}+\lambda_{j_1}}\cdots
p_{\lambda_{i_e}+\lambda_{j_e}}p_{\lambda_{i_{e+1}}}p_{\lambda_{i_{e+2}}}\cdots
p_{\lambda_{i_o}}.$$ We claim that $A$ appears as a term of the
lowest degree in the power sum expansion of $Q_{\lambda}$. For this
purpose, we need to determine those matchings $\pi$ of
$\{1,2,\ldots,2m\}$  in \eqref{q-expand}, for which the power sum
expansion of the product
$Q_{(\lambda_{\pi_1},\lambda_{\pi_2})}\cdots
Q_{(\lambda_{\pi_{2m-1}},\lambda_{\pi_{2m}})}$ contains $A$ as a
term of the lowest degree.

By Lemma \ref{lemma1}, if the $p_{\lambda_{i_1}+\lambda_{j_1}}$
appears as a factor in the lowest degree terms of the power sum
expansion of $Q_{(\lambda_{\pi_1},\lambda_{\pi_2})}\cdots
Q_{(\lambda_{\pi_{2m-1}},\lambda_{\pi_{2m}})}$, then we have
$\{\pi_1,\pi_2\}=\{i_1,j_1\}$. Iterating this argument, we see that
if  $p_{\lambda_{i_1}+\lambda_{j_1}}\cdots
p_{\lambda_{i_e}+\lambda_{j_e}}$ appears as a factor in the lowest
degree terms of $Q_{(\lambda_{\pi_1},\lambda_{\pi_2})}\cdots
Q_{(\lambda_{\pi_{2m-1}},\lambda_{\pi_{2m}})}$, then we have
$$\{\pi_1,\pi_2\}=\{i_1,j_1\},\ldots,\{\pi_{2e-1},\pi_{2e}\}=\{i_e,j_e\}.$$
It remains to determine the ordered pairs
$$\{(\pi_{2e+1},\pi_{2e+2}),\ldots,(\pi_{2m-1},\pi_{2m})\}.$$
 By the same argument as in Case (5) of the proof of Lemma \ref{len2}, for any
$e+1\leq k<l\leq o$, the term
$p_{\lambda_{i_{k}}}p_{\lambda_{i_{l}}}$ appears as a term of the
lowest degree in the power sum expansion of
$Q_{(\lambda_{i_k},\lambda_{i_l})}$.  Moreover, if the power sum
symmetric function $p_{\lambda_{i_{e+1}}}p_{\lambda_{i_{e+2}}}\cdots
p_{\lambda_{i_o}}$ appears as a term of the lowest degree in the
power sum expansion of the product
$Q_{(\lambda_{\pi_{2e+1}},\lambda_{\pi_{2e+2}})}\cdots
Q_{(\lambda_{\pi_{2m-1}},\lambda_{\pi_{2m}})}$, then the composition
of the pairs
$\{(\pi_{2e+1},\pi_{2e+2}),\ldots,(\pi_{2m-1},\pi_{2m})\}$ could be
any matching of $\{1,2,\ldots,2m\}/\{i_1,j_1,\ldots,i_e,j_e\}$.

To summarize, there are $(2(m-e)-1)!!$ matchings $\pi$ such that $A$
appears as a term of the lowest degree in the power sum expansion of
the product $Q_{(\lambda_{\pi_1},\lambda_{\pi_2})}\cdots
Q_{(\lambda_{\pi_{2m-1}},\lambda_{\pi_{2m}})}$. Combining Corollary
\ref{2odd} and Theorem \ref{conn}, we find that the coefficient of
$p_{\lambda_{i_k}}p_{\lambda_{i_l}}$ $(e+1\leq k<l\leq o)$ in the
power sum expansion of $Q_{(\lambda_{i_k}, \lambda_{i_l})}$ is
$-\frac{4}{\lambda_{i_k}\lambda_{i_l}}$. It follows that the
coefficient of $A$ in the expansion of the product
$Q_{(\lambda_{\pi_1},\lambda_{\pi_2})}\cdots
Q_{(\lambda_{\pi_{2m-1}}, \lambda_{\pi_{2m}})}$ is independent of
the choice of  $\pi$. Since $(2(m-e)-1)!!$ is an odd number, the
term $A$ will not vanish in the expansion of $Q_{\lambda}$. Note
that the degree of $A$ is $e+(o-e)=o,$ which is equal to ${\rm
srank}(\lambda)$, as desired.

Similarly, we consider the case $e>o$. In this case, we aim to show
that ${\rm srank}(\lambda)=e.$ By Theorem \ref{min bar}, if
$\lambda_{2m}>0$, i.e., $\ell(\lambda)=2m$, then we have
$${\rm srank}(\lambda)=\max(o,e+0)=e.$$ If $\lambda_{2m}=0$, i.e.,
$\ell(\lambda)=2m-1$, then we still have
$${\rm srank}(\lambda)=\max(o,(e-1)+1)=e.$$

Let
$$B=p_{\lambda_{i_1}+\lambda_{j_1}}\cdots
p_{\lambda_{i_o}+\lambda_{j_o}}p_{\lambda_{j_{o+1}}-1}p_{\lambda_{j_{o+2}}+1}\cdots
p_{\lambda_{j_{e-1}}-1}p_{\lambda_{j_e}+1}.$$ We proceed to prove
that $B$ appears as a term of the lowest degree in the power sum
expansion of $Q_{\lambda}$. Applying Lemma \ref{lemma1} repeatedly,
we deduce that if $p_{\lambda_{i_1}+\lambda_{j_1}}\cdots
p_{\lambda_{i_o}+\lambda_{j_o}}$ appears as a factor in the lowest
degree terms of the product
$Q_{(\lambda_{\pi_1},\lambda_{\pi_2})}\cdots
Q_{(\lambda_{\pi_{2m-1}},\lambda_{\pi_{2m}})}$, then
\begin{equation}\label{match1}
\{\pi_1,\pi_2\}=\{i_1,j_1\},\ldots,\{\pi_{2o-1},\pi_{2o}\}=\{i_o,j_o\}.
\end{equation}
On the other hand, iteration of  Lemma \ref{lemma2} reveals that if
the power sum symmetric function
$p_{\lambda_{j_{o+1}}-1}p_{\lambda_{j_{o+2}}+1}\cdots
p_{\lambda_{j_{e-1}}-1}p_{\lambda_{j_e}+1}$ appears as a term of the
lowest degree in the power sum expansion of
$Q_{(\lambda_{\pi_{2o+1}},\lambda_{\pi_{2o+2}})}\cdots
Q_{(\lambda_{\pi_{2m-1}},\lambda_{\pi_{2m}})}$, then
\begin{equation}\label{match2}
\{\pi_{2o+1},\pi_{2o+2}\}=\{j_{o+1},j_{o+2}\},\ldots,\{\pi_{2m-1},\pi_{2m}\}=\{j_{e-1},j_e\}.
\end{equation}
Therefore, if $B$ appears as a term of the lowest degree in the
power sum expansion of $Q_{(\lambda_{\pi_1},\lambda_{\pi_2})}\cdots
Q_{(\lambda_{\pi_{2m-1}},\lambda_{\pi_{2m}})}$, then the matching
$\pi$ is uniquely determined by \eqref{match1} and \eqref{match2}.
Note that the degree of $B$ is $e$, which coincides with ${\rm
srank}(\lambda)$.

Since there is always a term of degree ${\rm srank}(\lambda)$ in the
power sum expansion of $Q_\lambda$, the theorem follows. \qed

\section{Skew Schur's $Q$-functions}

In this section, we show that the srank ${\rm srank}(\lambda/\mu)$
is a lower bound of the lowest degree of the terms in the power sum
expansion of the skew Schur's $Q$-function $Q_{\lambda/\mu}$. Note
that Clifford's conjecture does not hold for skew shapes.

We first recall a definition of the skew Schur's $Q$-function in
terms of strip tableaux. The concept of strip tableaux were
introuduced by Stembridge \cite{stemb1988} to describe the Morris
rule for the evaluation of irreducible spin characters. Given a skew
partition $\lambda/\mu$, the \emph{$j$-th diagonal} of the skew
shifted diagram $S(\lambda/\mu)$ is defined as the set of squares
$(1,j), (2, j+1), (3, j+2), \ldots$ in $S(\lambda/\mu)$. A skew
diagram $S(\lambda/\mu)$ is called a \emph{strip} if it is rookwise
connected and each diagonal contains at most one box. The
\emph{height} $h$ of a strip is defined to be the number of rows it
occupies. A \emph{double strip} is a skew diagram formed by the
union of two strips which both start on the diagonal consisting of
squares $(j,j)$. The \emph{depth} of a double strip is defined to be
$\alpha+\beta$ if it has $\alpha$ diagonals of length two and its
diagonals of length one occupy $\beta$ rows. A \emph{strip tableau}
of shape $\lambda/\mu$ and type $\pi=(\pi_1,\ldots,\pi_k)$ is
defined to be a sequence of shifted diagrams
$$S(\mu)=S(\lambda^0)\subseteq S(\lambda^1)\subseteq \cdots \subseteq S(\lambda^k)=S(\lambda)$$
with $|\lambda^i/\lambda^{i-1}|=\pi_i$ ($1\leq i\leq k$) such that
each skew shifted diagram $S(\lambda^i/\lambda^{i-1})$ is either a
strip or a double strip.

The skew Schur's $Q$-function can be defined as the weight
generating function of strip tableaux in the following way. For a
strip of height $h$ we assign the weight $(-1)^{h-1}$, and for a
double strip of depth $d$ we assign the weight $2(-1)^{d-1}$. The
weight of a strip tableau $T$, denoted $wt(T)$, is the product of
the weights of strips and double strips of which $T$ is composed.
Then the skew Schur's $Q$-function $Q_{\lambda/\mu}$ is given by
\begin{equation}
Q_{\lambda/\mu}=\sum_{\pi\in \mathcal{P}^o(|\lambda/\mu|)}\sum_{T}
2^{\ell(\pi)}wt(T)\frac{p_{\pi}}{z_{\pi}},
\end{equation}
where $T$ ranges over all strip tableaux $T$ of shape $\lambda/\mu$
and type $\pi$, see \cite[Theorem 5.1]{stemb1988}.

J$\rm{\acute{o}}$zefiak and Pragacz \cite{jozpra1991} obtained the
following Pfaffian formula for the skew Schur's $Q$-function.
\begin{thm}[\cite{jozpra1991}]\label{skewpf}
Let $\lambda, \mu$ be strict partitions with $m=\ell(\lambda)$,
$n=\ell(\mu)$, $\mu\subset \lambda$, and let $M(\lambda,\mu)$ denote
the skew-symmetric matrix
$$\begin{pmatrix}
A & B\\ -B^t & 0
\end{pmatrix},
$$
where $A=(Q_{(\lambda_i,\lambda_j)})$ and
$B=(Q_{(\lambda_i-\mu_{n+1-j})})$.

Then
\begin{itemize}
\item[(1)] if $m+n$ is even, we have $Q_{\lambda/\mu}={\rm
Pf}(M(\lambda,\mu))$;

\item[(2)] if $m+n$ is odd, we have $Q_{\lambda/\mu}={\rm Pf}(M(\lambda,\mu^\prime))$, where $\mu^\prime=(\mu_1,\cdots,\mu_n,
0)$.
\end{itemize}
\end{thm}
A combinatorial proof of the above theorem was given by Stembridge
\cite{stemb1990} in terms of lattice paths, and later, Hamel
\cite{hamel1996} gave an interesting generalization by using the
border strip decompositions of the shifted diagram.

Given a skew partition $\lambda/\mu$, Clifford \cite{cliff2003}
constructed a bijection between skew bar tableaux of shape
$\lambda/\mu$ and skew strip tableaux of the same shape, which
preserves the type of the tableau. Using this bijection, it is
straightforward to derive the following result.
\begin{prop}
The terms of the lowest degree in $Q_{\lambda/\mu}$ have degree at
least ${\rm srank}(\lambda/\mu)$.
\end{prop}

Different from the case of non-skew shapes, in general, the lowest
degree terms in $Q_{\lambda/\mu}$ do not have the degree  ${\rm
srank}(\lambda/\mu)$. For example, take the skew partition
$(4,3)/3$. It is easy to see that ${\rm srank}((4,3)/3)=2$. While,
using Theorem \ref{skewpf} and Stembridge's SF Package for Maple
\cite{stem2}, we obtain that
\begin{equation}
Q_{(4,3)/3}={\rm Pf}
\begin{pmatrix}0 & Q_{(4,3)} & Q_{(4)} & Q_{(1)}\\[5pt]
Q_{(3,4)} & 0 & Q_{(3)} & Q_{(0)}\\[5pt]
-Q_{(4)} & -Q_{(3)}& 0  & 0\\[5pt]
-Q_{(1)} & -Q_{(0)}& 0  & 0
\end{pmatrix}=2p_1^4.
\end{equation}
This shows that the lowest degree of $Q_{(4,3)/3}$ equals 4, which
is strictly greater than ${\rm srank}((4,3)/3)$.

\section{The srank of skew partitions} \label{sect4}

In this section, we  present an algorithm to determine the srank for
the skew partition $\lambda/\mu$. In fact, the algorithm leads to a
configuration of $0$'s. To obtain the srank of a skew partition, we
need to minimize the number of bars by adjusting the positions of
$0$'s. Given a configuration $\mathcal{C}$ of $0$'s in the shifted
diagram $S(\lambda)$, let
$$\kappa(\mathcal{C})=o_s+2e_s+\max(o_r,e_r+((e_r+o_r)\ \mathrm{mod}\  2)),$$
where $o_ r$ (resp. $e_r$) counts the number of nonempty rows in
which there are an odd (resp. even) number of squares  and no
squares are filled with $0$, and $o_s$ (resp. $e_s$) records the
number of rows in which at least one square is filled with $0$ but
there are an odd (resp. nonzero even) number of blank squares.

If there exists at least one bar tableau of type $\lambda/\mu$ under
some configuration $\mathcal{C}$, we say that $\mathcal{C}$ is
\emph{admissible}. For a fixed configuration $\mathcal{C}$, each row
is one of the following eight possible types:
\begin{itemize}
\item[(1)] an even row bounded by an even number of $0$'s, denoted
$(e,e)$,

\item[(2)] an odd row bounded by an even number of $0$'s, denoted
$(e,o)$,

\item[(3)] an odd row bounded by an odd number of $0$'s, denoted
$(o,e)$,

\item[(4)] an even row bounded by an odd number of $0$'s, denoted
$(o,o)$,

\item[(5)] an even row without $0$'s, denoted $(\emptyset,e)$,

\item[(6)] an  odd row without $0$'s, denoted $(\emptyset,o)$,

\item[(7)] an even row filled with $0$'s, denoted $(e,
\emptyset)$,

\item[(8)] an  odd row filled with $0$'s, denoted $(o,
\emptyset)$.

\end{itemize}

Given two rows with respective types $s$ and $s'$ for some
configuration $\mathcal{C}$, if we can obtain a new configuration
$\mathcal{C}'$ by exchanging the locations of $0$'s in these two
rows such that their new types are $t$ and $t'$ respectively, then
denote it by $\mathcal{C}'=\mathcal{C}\left( \left[\tiny{{s \atop
s'}}\right] \rightarrow \left[\tiny{{{t} \atop
{t'}}}\right]\right)$. Let $o_r,e_r,o_s,e_s$ be defined as above
corresponding to configuration $\mathcal{C}$, and let
$o_r',e_r',o_s',e_s'$ be those of $\mathcal{C}'$.

In the following we will show that how the quantity
$\kappa(\mathcal{C})$ changes when exchanging the locations of $0$'s
in $\mathcal{C}$.

\begin{lem}\label{varyzero1-1}
If $\mathcal{C}'=\mathcal{C}\left( \left[\tiny{{s \atop s'}}\right]
\rightarrow \left[\tiny{{s \atop s'}}\right]\right)$ or
$\mathcal{C}'=\mathcal{C}\left( \left[\tiny{{s \atop s'}}\right]
\rightarrow \left[\tiny{{{s'} \atop s}}\right]\right)$, i.e., the
types of the two involved rows are remained or exchanged, where
$s,s'$ are any two possible types, then $\kappa({\mathcal{C}'})=
\kappa({\mathcal{C}})$.
\end{lem}

\begin{lem}\label{varyzero1-6} If
$\mathcal{C}'=\mathcal{C}\left( \left[\tiny{{{(e,e)}
\atop\rule{0pt}{10pt} {(\emptyset,o)}}}\right] \rightarrow
\left[\tiny{{{(\emptyset,e)} \atop\rule{0pt}{10pt}
{(e,o)}}}\right]\right)$, then $\kappa({\mathcal{C}'})\leq
\kappa({\mathcal{C}})$.
\end{lem}

\pf In this case we have
$$o_s^\prime=o_s+1, \quad e_s^\prime=e_s-1, \quad
o_r^\prime=o_r-1,\quad e_r^\prime=e_r+1.$$ Note that
$o_r+e_r=\ell(\lambda)-\ell(\mu)$. Now there are two cases to
consider.

\noindent \textbf{Case I.} The skew partition $\lambda/\mu$
satisfies that $\ell(\lambda)-\ell(\mu) \equiv 0\ (\mathrm{mod}\
2)$.
\begin{itemize}
\item[(1)] If $o_r\leq e_r$, then $o_r^\prime\leq e_r^\prime$ and
\begin{align*}
\kappa(\mathcal{C})&=o_s+2e_s+e_r,\\
\kappa(\mathcal{C}')&=o_s+1+2(e_s-1)+e_r^\prime=o_s+2e_s+e_r=\kappa(\mathcal{C}).
\end{align*}

\item[(2)] If $o_r\geq e_r+2$, then $o_r^\prime=o_r-1\geq
e_r+1=e_r^\prime$ and
\begin{align*}
\kappa(\mathcal{C})&=o_s+2e_s+o_r,\\
\kappa(\mathcal{C}')&=o_s+2e_s-1+o_r^\prime=o_s+2e_s+o_r-2<\kappa(\mathcal{C}).
\end{align*}
\end{itemize}

\noindent \textbf{Case II.} The skew partition $\lambda/\mu$
satisfies that $\ell(\lambda)-\ell(\mu) \equiv 1\ (\mathrm{mod}\
2)$.
\begin{itemize}
\item[(1)] If $o_r\leq e_r+1$, then $o_r^\prime<e_r^\prime$ and
\begin{align*}
\kappa(\mathcal{C})&=o_s+2e_s+e_r+1,\\
\kappa(\mathcal{C}')&=o_s+2e_s-1+e_r^\prime+1=o_s+2e_s+e_r+1=\kappa(\mathcal{C}).
\end{align*}

\item[(2)] If $o_r\geq e_r+3$, then $o_r^\prime=o_r-1\geq
e_r+2>e_r^\prime$ and
\begin{align*}
\kappa(\mathcal{C})&=o_s+2e_s+o_r,\\
\kappa(\mathcal{C}')&=o_s+2e_s-1+o_r^\prime=o_s+2e_s+o_r-2<\kappa(\mathcal{C}).
\end{align*}
\end{itemize}
Therefore, the inequality $\kappa({\mathcal{C}'})\leq
\kappa({\mathcal{C}})$ holds under the assumption. \qed

\begin{lem}\label{varyzero2-6} If
$\mathcal{C}'=\mathcal{C}\left( \left[\tiny{{{(o,e)}
\atop\rule{0pt}{10pt} {(\emptyset,e)}}}\right] \rightarrow
\left[\tiny{{{(\emptyset,o)} \atop\rule{0pt}{10pt}
{(o,o)}}}\right]\right)$, then $\kappa({\mathcal{C}'})\leq
\kappa({\mathcal{C}})$.
\end{lem}

\pf In this case we have
$$o_s^\prime=o_s+1, \quad e_s^\prime=e_s-1, \quad
o_r^\prime=o_r+1,\quad e_r^\prime=e_r-1.$$ Now there are two
possibilities.

\noindent \textbf{Case I.} The skew partition $\lambda/\mu$
satisfies that $\ell(\lambda)-\ell(\mu) \equiv 0\ (\mathrm{mod}\
2)$.
\begin{itemize}
\item[(1)] If $o_r\leq e_r-2$, then $o_r^\prime\leq e_r^\prime$
and
\begin{align*}
\kappa(\mathcal{C})&=o_s+2e_s+e_r,\\
\kappa(\mathcal{C}')&=o_s+1+2(e_s-1)+e_r^\prime=o_s+2e_s+e_r-2<\kappa(\mathcal{C}).
\end{align*}

\item[(2)] If $o_r\geq e_r$, then $o_r^\prime=o_r+1>
e_r-1=e_r^\prime$ and
\begin{align*}
\kappa(\mathcal{C})&=o_s+2e_s+o_r,\\
\kappa(\mathcal{C}')&=o_s+2e_s-1+o_r^\prime=o_s+2e_s+o_r=\kappa(\mathcal{C}).
\end{align*}
\end{itemize}

\noindent \textbf{Case II.} The skew partition $\lambda/\mu$
satisfies that $\ell(\lambda)-\ell(\mu) \equiv 1\ (\mathrm{mod}\
2)$.
\begin{itemize}
\item[(1)] If $o_r\leq e_r-3$, then $o_r^\prime<e_r^\prime$ and
\begin{align*}
\kappa(\mathcal{C})&=o_s+2e_s+e_r+1,\\
\kappa(\mathcal{C}')&=o_s+2e_s-1+e_r^\prime+1=o_s+2e_s+e_r-1<\kappa(\mathcal{C}).
\end{align*}

\item[(2)] If $o_r\geq e_r-1$, then $o_r^\prime=o_r+1>
e_r-1=e_r^\prime$ and
\begin{align*}
\kappa(\mathcal{C})&=o_s+2e_s+o_r,\\
\kappa(\mathcal{C}')&=o_s+2e_s-1+o_r^\prime=o_s+2e_s+o_r=\kappa(\mathcal{C}).
\end{align*}
\end{itemize}
In both  cases we have $\kappa({\mathcal{C}'})\leq
\kappa({\mathcal{C}})$, as required. \qed

\begin{lem}\label{varyzero1-4} If
$\mathcal{C}'=\mathcal{C}\left( \left[\tiny{{{(e,e)}
\atop\rule{0pt}{10pt} {(o,e)}}}\right] \rightarrow
\left[\tiny{{{(o,o)} \atop\rule{0pt}{10pt}
{(e,o)}}}\right]\right)$, then $\kappa({\mathcal{C}'})<
\kappa({\mathcal{C}})$.
\end{lem}

\pf In this case, we have
$$o_s^\prime=o_s+2, \quad e_s^\prime=e_s-2, \quad
o_r^\prime=o_r,\quad e_r^\prime=e_r.$$ Therefore,
$$
\kappa(\mathcal{C}')=o_s'+2e_s'+\max(o_r',e_r'+((e_r'+o_r')\
\mathrm{mod}\ 2))=\kappa(\mathcal{C})-2.
$$
The desired inequality immediately follows.  \qed


\begin{lem}\label{varyzero3-5}
If $\mathcal{C}'=\mathcal{C}\left( \left[\tiny{{{(e,o)}
\atop\rule{0pt}{10pt} {(\emptyset,e)}}}\right] \rightarrow
\left[\tiny{{{(\emptyset,o)} \atop\rule{0pt}{10pt}
{(e,\emptyset)}}}\right]\right)$, then $\kappa({\mathcal{C}'})\leq
\kappa({\mathcal{C}})$.
\end{lem}

\pf Under this transformation we have
$$o_s^\prime=o_s-1, \quad e_s^\prime=e_s, \quad
o_r^\prime=o_r+1,\quad e_r^\prime=e_r-1.$$ Since
$o_r+e_r=\ell(\lambda)-\ell(\mu)$ is invariant, there are two cases.

\noindent \textbf{Case I.} The skew partition $\lambda/\mu$
satisfies that $\ell(\lambda)-\ell(\mu) \equiv 0\ (\mathrm{mod}\
2)$.
\begin{itemize}

\item[(1)] If $o_r\geq e_r$, then $o_r^\prime\geq e_r^\prime$ and
\begin{align*}
\kappa(\mathcal{C}')=o_s^\prime+2e_s^\prime+o_r^\prime
=o_s-1+2e_s+o_r+1 =\kappa(\mathcal{C}).
\end{align*}

\item[(2)] If $o_r\leq e_r-2$, then $o_r^\prime=o_r+1\leq
e_r-1=e_r^\prime$ and
\begin{align*}
\kappa(\mathcal{C}')=o_s-1+2e_s+e_r^\prime=o_s+2e_s+e_r-2<\kappa(\mathcal{C}).
\end{align*}

\end{itemize}

\noindent \textbf{Case II.} The skew partition $\lambda/\mu$
satisfies that $\ell(\lambda)-\ell(\mu) \equiv 1\ (\mathrm{mod}\
2)$.
\begin{itemize}
\item[(1)] If $o_r\geq e_r+1$, then $o_r^\prime=o_r+1\geq
e_r+2>e_r^\prime+1$ and
\begin{align*}
\kappa(\mathcal{C}')=o_s^\prime+2e_s^\prime+o_r^\prime=o_s+2e_s+o_r=\kappa(\mathcal{C}).
\end{align*}

\item[(2)] If $o_r\leq e_r-1$, then $o_r^\prime=o_r+1\leq
e_r=e_r^\prime+1$ and
\begin{align*}
\kappa(\mathcal{C}')=o_s^\prime+2e_s^\prime+e_r^\prime+1=o_s+2e_s+e_r-1<\kappa(\mathcal{C}).
\end{align*}
\end{itemize}
Hence the proof is complete. \qed


\begin{lem}\label{varyzero5-8}
If $\mathcal{C}'=\mathcal{C}\left( \left[\tiny{{{(o,o)}
\atop\rule{0pt}{10pt} {(\emptyset,o)}}}\right] \rightarrow
\left[\tiny{{{(\emptyset,e)} \atop\rule{0pt}{10pt}
{(o,\emptyset)}}}\right]\right)$, then $\kappa({\mathcal{C}'})\leq
\kappa({\mathcal{C}})$.
\end{lem}

\pf In this case we have
$$o_s^\prime=o_s-1, \quad e_s^\prime=e_s, \quad
o_r^\prime=o_r-1,\quad e_r^\prime=e_r+1.$$ There are two
possibilities:

 \noindent \textbf{Case I.} The skew
partition $\lambda/\mu$ satisfies that $\ell(\lambda)-\ell(\mu)
\equiv 0\ (\mathrm{mod}\ 2)$.
\begin{itemize}

\item[(1)] If $o_r\geq e_r+2$, then $o_r^\prime=o_r-1\geq
e_r+1=e_r^\prime$ and
\begin{align*}
\kappa(\mathcal{C}')=o_s^\prime+2e_s^\prime+o_r^\prime
=o_s-1+2e_s+o_r-1<\kappa(\mathcal{C}).
\end{align*}

\item[(2)] If $o_r\leq e_r$, then $o_r^\prime=o_r-1\leq
e_r-1<e_r^\prime$ and
\begin{align*}
\kappa(\mathcal{C}')=o_s-1+2e_s+e_r^\prime=o_s+2e_s+e_r=\kappa(\mathcal{C}).
\end{align*}

\end{itemize}

\noindent \textbf{Case II.} The skew partition $\lambda/\mu$
satisfies that $\ell(\lambda)-\ell(\mu) \equiv 1\ (\mathrm{mod}\
2)$.
\begin{itemize}
\item[(1)] If $o_r\geq e_r+3$, then $o_r^\prime=o_r-1\geq
e_r+2=e_r^\prime+1$ and
\begin{align*}
\kappa(\mathcal{C}')=o_s^\prime+2e_s^\prime+o_r^\prime=o_s+2e_s+o_r-2<
\kappa(\mathcal{C}).
\end{align*}

\item[(2)] If $o_r\leq e_r+1$, then $o_r^\prime=o_r-1\leq
e_r<e_r^\prime+1$ and
\begin{align*}
\kappa(\mathcal{C}')=o_s^\prime+2e_s^\prime+e_r^\prime+1=o_s+2e_s+e_r+1=\kappa(\mathcal{C}).
\end{align*}
\end{itemize}
Therefore, in both cases we have $\kappa({\mathcal{C}'})\leq
\kappa({\mathcal{C}})$. \qed

\begin{lem}\label{varyzero2-3}
If $\mathcal{C}'=\mathcal{C}\left( \left[\tiny{{{(e,o)}
\atop\rule{0pt}{10pt} {(o,o)}}}\right] \rightarrow
\left[\tiny{{{(o,e)} \atop\rule{0pt}{10pt}
{(e,\emptyset)}}}\right]\right)$ or $\mathcal{C}'=\mathcal{C}\left(
\left[\tiny{{{(o,o)} \atop\rule{0pt}{10pt}  {(e,o)}}}\right]
\rightarrow \left[\tiny{{{(e,e)} \atop\rule{0pt}{10pt}
{(o,\emptyset)}}}\right]\right)$, then $\kappa({\mathcal{C}'})=
\kappa({\mathcal{C}})$.
\end{lem}

\pf In each case we have
$$o_s^\prime=o_s-2, \quad e_s^\prime=e_s+1, \quad
o_r^\prime=o_r,\quad e_r^\prime=e_r.$$ Therefore
$$
\kappa(\mathcal{C}')=o_s'+2e_s'+\max(o_r',e_r'+((e_r'+o_r')\
\mathrm{mod}\ 2))=\kappa(\mathcal{C}),
$$
as desired.  \qed


\begin{lem}\label{varyzero1-7}
If $\mathcal{C}'$ is one of the following possible cases:
\[
\begin{array}{ccc}
\mathcal{C}\left( \left[\tiny{{{(e,e)} \atop \rule{0pt}{10pt}
{(e,e)}}}\right] \rightarrow \left[\tiny{{{(e,e)} \atop
\rule{0pt}{10pt} {(e,\emptyset)}}}\right]\right), &
\mathcal{C}\left( \left[\tiny{{{(e,e)} \atop \rule{0pt}{10pt}
{(o,o)}}}\right] \rightarrow \left[\tiny{{{(o,o)}
\atop\rule{0pt}{10pt} {(e,\emptyset)}}}\right]\right), &
\mathcal{C}\left( \left[\tiny{{{(e,o)} \atop \rule{0pt}{10pt}
{(e,e)}}}\right] \rightarrow \left[\tiny{{{(e,o)} \atop
\rule{0pt}{10pt} {(e,\emptyset)}}}\right]\right),\\[10pt]
\mathcal{C}\left( \left[\tiny{{{(e,e)} \atop\rule{0pt}{10pt}
{(\emptyset,e)}}}\right] \rightarrow \left[\tiny{{{(\emptyset,e)}
\atop\rule{0pt}{10pt} {(e,\emptyset)}}}\right]\right), &
\mathcal{C}\left( \left[\tiny{{{(o,o)} \atop \rule{0pt}{10pt}
{(o,e)}}}\right] \rightarrow \left[\tiny{{{(o,o)} \atop
\rule{0pt}{10pt}  {(o,\emptyset)}}}\right]\right), &
\mathcal{C}\left( \left[\tiny{{ {(o,e)} \atop\rule{0pt}{10pt}
{(e,o)}}}\right] \rightarrow \left[\tiny{{{(e,o)}
\atop\rule{0pt}{10pt} {(o,\emptyset)}}}\right]\right),\\[10pt]
\mathcal{C}\left( \left[\tiny{{{(o,e)} \atop\rule{0pt}{10pt}
{(o,e)}}}\right] \rightarrow \left[\tiny{{{(o,e)}
\atop\rule{0pt}{10pt} {(o,\emptyset)}}}\right]\right), &
\mathcal{C}\left( \left[\tiny{{{(o,e)} \atop\rule{0pt}{10pt}
{(\emptyset,o)}}}\right] \rightarrow \left[\tiny{{{(\emptyset,o)}
\atop \rule{0pt}{10pt} {(o,\emptyset)}}}\right]\right),&
\end{array}
\]
 then
$\kappa({\mathcal{C}'})< \kappa({\mathcal{C}})$.
\end{lem}

\pf In each case we have
$$o_s^\prime=o_s, \quad e_s^\prime=e_s-1, \quad
o_r^\prime=o_r,\quad e_r^\prime=e_r.$$ Therefore
$$
\kappa(\mathcal{C}')=o_s'+2e_s'+\max(o_r',e_r'+((e_r'+o_r')\
\mathrm{mod}\ 2))<\kappa(\mathcal{C}),
$$
as required. \qed

Note that Lemmas \ref{varyzero1-1}-\ref{varyzero1-7} cover all
possible transformations of exchanging the locations of $0$'s in two
involved rows. Lemmas \ref{varyzero1-6}-\ref{varyzero1-4} imply
that, to minimize the number of bars, we should put $0$'s in the
skew shifted diagram such that there are as more as possible rows
for which the first several squares are filled with $0$'s and then
followed by an odd number of blank squares. Meanwhile, from Lemmas
\ref{varyzero3-5}-\ref{varyzero1-7} we know that the number of rows
fully filled with $0$'s should be as more as possible. Based on
these observations, we have the following algorithm to determine the
location of $0$'s for a given skew partition $\lambda/\mu$, where
both $\lambda$ and $\mu$ are strict partitions. Using this algorithm
we will obtain a shifted diagram with some squares filled with $0$'s
such that the corresponding quantity $\kappa(\mathcal{C})$ is
minimized.  This property allows us to determine the srank of
$\lambda/\mu$.

\noindent{\textbf{The Algorithm for Determining the Locations of
$0$'s:}}
\begin{itemize}
\item[(S1)] Let $\mathcal{C}_1=S(\lambda)$ be the
initial configuration of $\lambda/\mu$ with blank square. Set $i=1$
and $J=\{1,\ldots,\ell(\lambda)\}$.

\item[(S2)] For $i\leq \ell(\mu)$, iterate the following procedure:
    \begin{itemize}
     \item[(A)] If $\mu_i=\lambda_j$ for some $j\in J$, then we fill the
$j$-th row of $\mathcal{C}_i$ with $0$.

     \item[(B)] If $\mu_i\neq \lambda_j$ for any $j\in J$, then
     there are two possibilities.
        \begin{itemize}
         \item[(B1)] $\lambda_j-\mu_i$ is odd for some $j\in J$
         and $\lambda_j>\mu_i$. Then we take the largest such $j$ and
         fill the leftmost $\mu_i$ squares with $0$ in the $j$-th row of
         $\mathcal{C}_i$.

         \item[(B2)] $\lambda_j-\mu_i$ is even for any $j\in J$
         and $\lambda_j>\mu_i$. Then we take the largest such $j$ and
         fill the leftmost $\mu_i$ squares by $0$ in the $j$-th row of
         $\mathcal{C}_i$.
         \end{itemize}
    \end{itemize}

    Denote the new configuration
by $\mathcal{C}_{i+1}$. Set $J=J\backslash \{j\}$.

\item[(S3)] Set $\mathcal{C}^{*}=\mathcal{C}_{i}$, and we get the
desired configuration.
\end{itemize}

It should be emphasized  that although the above algorithm does not
necessarily generate a bar tableau, it is sufficient for the
computation of the srank of a skew partition.

Using the arguments in the proofs of Lemmas
\ref{varyzero1-1}-\ref{varyzero1-7}, we can derive the following
crucial property of the configuration $\mathcal{C}^*$. The proof is
omitted since it is tedious and  straightforward.

\begin{prop}\label{prop-min}
For any configuration ${\mathcal{C}}$ of $0$'s in the skew shifted
diagram of $\lambda/\mu$, we have $\kappa({\mathcal{C}^*})\leq
\kappa({\mathcal{C}})$.
\end{prop}


\begin{thm}\label{number of skew}
Given a skew partition $\lambda/\mu$, let $\mathcal{C}^*$ be the
configuration of $0$'s obtained by applying the algorithm described
above. Then
\begin{equation}\label{srank}
{\rm srank}(\lambda/\mu)=\kappa({\mathcal{C}^*}).
\end{equation}
\end{thm}
\pf Suppose that for the configuration ${\mathcal{C}^*}$ there are
$o_r^*$ rows of odd size with blank squares, and there are $o_s^*$
rows with at least one square filled with $0$ and an odd number of
squares filled with positive integers. Likewise we let $e_r^*$ and
$e_s^*$ denote the number of remaining rows. Therefore,
$$\kappa(\mathcal{C}^*)=o_s^*+2e_s^*+\max(o_r^*,e_r^*+((e_r^*+o_r^*)\ \mathrm{mod}\  2)).$$
Since for each configuration $\mathcal{C}$ the number of bars in a
minimal bar tableau is greater than or equal to
$\kappa({\mathcal{C}})$, by Proposition \ref{prop-min}, it suffices
to confirm the existence of a skew bar tableau, say $T$, with
$\kappa({\mathcal{C}^*})$ bars.

Note that it is possible that the configuration ${\mathcal{C}^*}$ is
not admissible. The key idea of our proof is to move $0$'s in the
diagram such that the resulting configuration ${\mathcal{C}'}$ is
admissible and $\kappa({\mathcal{C}'})=\kappa({\mathcal{C}^*})$. To
achieve this goal, we will use the numbers
$\{1,2,\ldots,\kappa({\mathcal{C}^*})\}$ to fill up the blank
squares of $\mathcal{C}^*$ guided by the rule that the bars of Type
$2$ or Type $3$ will occur  before bars of Type $1$.

Let us consider the rows without $0$'s, and there are two
possibilities:  (A) $o_r^*\geq e_r^*$,  (B) $o_r^*<e_r^*$.

In Case (A) we choose a row of even size and a row of odd size, and
fill up these two rows with $\kappa({\mathcal{C}^*})$ to generate a
bar of Type $3$. Then we continue to choose a row of even size and a
row of odd size, and fill up these two rows with
$\kappa({\mathcal{C}^*})-1$. Repeat this procedure until all even
rows are filled up. Finally, we fill the remaining rows of odd size
with $\kappa({\mathcal{C}^*})-e_r^*,
\kappa({\mathcal{C}^*})-e_r^*-1, \ldots,
\kappa({\mathcal{C}^*})-o_r^*+1$ to generate bars of Type $2$.

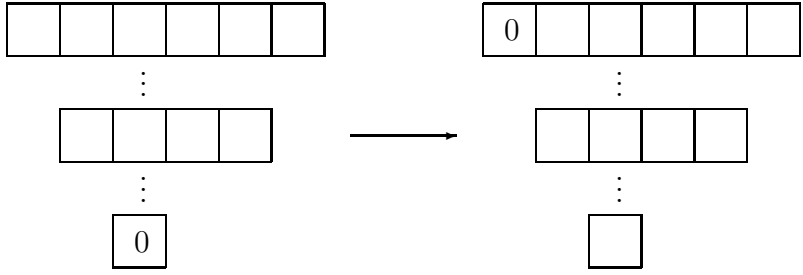
\begin{figure}[h,t]
\setlength{\unitlength}{1pt}
\begin{center}
\begin{picture}(300,100)
\put(40,0){\line(1,0){20}}\put(40,20){\line(1,0){20}}
\put(40,0){\line(0,1){20}}\put(60,0){\line(0,1){20}}
\put(50,25){$\vdots$}\put(20,40){\line(1,0){80}}
\put(20,60){\line(1,0){80}}
\multiput(20,40)(20,0){5}{\line(0,1){20}} \put(50,65){$\vdots$}
\multiput(0,80)(0,20){2}{\line(1,0){120}}
\multiput(0,80)(20,0){7}{\line(0,1){20}}\put(48,6){$0$}
\put(130,50){\vector(1,0){40}}
\put(220,0){\line(1,0){20}}\put(220,20){\line(1,0){20}}
\put(220,0){\line(0,1){20}}\put(240,0){\line(0,1){20}}
\put(230,25){$\vdots$}\put(200,40){\line(1,0){80}}
\put(200,60){\line(1,0){80}}
\multiput(200,40)(20,0){5}{\line(0,1){20}} \put(230,65){$\vdots$}
\multiput(180,80)(0,20){2}{\line(1,0){120}}
\multiput(180,80)(20,0){7}{\line(0,1){20}}\put(188,86){$0$}
\end{picture}
\end{center}
\caption{Vacating the unique square at the bottom of the
diagram}\label{case2-2}
\end{figure}

In Case (B) we choose the row with the $i$-th smallest even size and
the row with the $i$-th smallest  odd size and fill their squares
with the  number $\kappa({\mathcal{C}^*})-i+1$ for
$i=1,\ldots,o_r^*$. In this way, we obtain $o_r^*$ bars of Type $3$.
Now consider the remaining rows of even size without $0$'s. There
are two subcases.
\begin{itemize}
\item[(B1)] The remaining diagram, obtained by removing the
previous $o_r^*$ bars of Type $3$, does not contain any row with
only one square. Under this assumption, it is possible to fill the
squares of a row of even size with the number
$\kappa({\mathcal{C}^*})-o_r^*$ except the leftmost square. This
operation will result in a bar of Type $1$. After removing this bar
from the diagram, we may combine this leftmost square of the current
row and another row of even size, if it exists, and to generate a
bar of Type $3$. Repeating this procedure until there are no more
rows of even size, we  obtain a sequence of bars of Type $1$ and
Type $3$. Evidently, there is a bar of Type $2$ with only one
square. To summarize, we have $\max(o_r^*,e_r^*+((e_r^*+o_r^*)\
\mathrm{mod}\ 2))$ bars.

\item[(B2)] The remaining diagram
contains a row composed of the unique square filled with $0$. In
this case, we will move this $0$ into the leftmost square of a row
of even size, see Figure \ref{case2-2}. Denote this new
configuration by $\mathcal{C}^{\prime}$, and  from Lemma
\ref{varyzero5-8} we see that
$\kappa({\mathcal{C}^*})=\kappa({\mathcal{C}^{\prime}})$. If we
start with ${\mathcal{C}'}$ instead of ${\mathcal{C}^*}$, by a
similar construction, we get $\max(o_r',e_r'+((e_r'+o_r')\
\mathrm{mod}\ 2))$ bars, occupying the rows without $0$'s in the
diagram.
\end{itemize}

Without loss of generality, we may assume that for the configuration
${\mathcal{C}^*}$ the rows without $0$'s in the diagram have been
occupied by the bars with the first
$\max(o_r^*,e_r^*+((e_r^*+o_r^*)\ \mathrm{mod}\ 2))$ positive
integers in the decreasing order, namely, $(\kappa({\mathcal{C}^*}),
\ldots, 2, 1, 0)$. By removing these bars and reordering the
remaining rows, we may get  a shifted diagram with which we can
continue the above procedure  to construct a bar tableau.

\begin{figure}[h,t]
\setlength{\unitlength}{1pt}
\begin{center}
\begin{picture}(340,40)
\put(20,0){\line(1,0){120}}
\multiput(0,20)(0,20){2}{\line(1,0){140}}
\multiput(20,0)(20,0){7}{\line(0,1){40}} \put(0,20){\line(0,1){20}}
\multiput(28,6)(20,0){6}{$0$}\multiput(8,26)(20,0){3}{$0$}
\put(150,20){\vector(1,0){40}} \put(220,0){\line(1,0){120}}
\multiput(200,20)(0,20){2}{\line(1,0){140}}
\multiput(220,0)(20,0){7}{\line(0,1){40}}
\put(200,20){\line(0,1){20}}
\multiput(228,6)(20,0){3}{$0$}\multiput(208,26)(20,0){6}{$0$}
\end{picture}
\end{center}
\caption{Interchanging the location of $0$'s in two neighbored
rows}\label{case2'}
\end{figure}
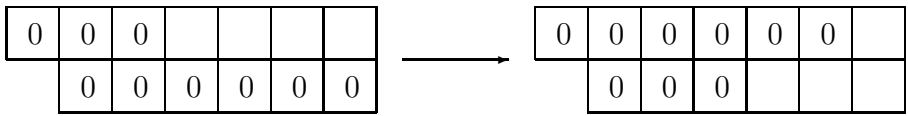

At this point, it is necessary to show that it is possible to use
$o_s^*+2e_s^*$ bars to fill this diagram. In doing so, we process
the rows from bottom to top. If the bottom row has an odd number of
blank squares, then we simply assign the symbol $o_s^*+2e_s^*$ to
these squares to produce a bar of Type $1$. If the bottom row are
completely filled with $0$'s, then we continue to deal with the row
above the bottom row. Otherwise, we fill the rightmost square of the
bottom row with $o_s^*+2e_s^*$ and the remaining squares with
$o_s^*+2e_s^*-1$. Suppose that we have filled $i$ rows from the
bottom and  all the involved bars have been removed from the
diagram. Then we consider the $(i+1)$-th row from the bottom. Let
$t$ denote the largest number not greater than $o_s^*+2e_s^*$ which
has not been used before. If all squares in the $(i+1)$-th row are
filled with $0$'s, then we continue to deal with the $(i+2)$-th row.
If the number of blank squares in the $(i+1)$-th row is odd, then we
fill these squares with $t$. If the number of blank squares in the
$(i+1)$-th row is even, then we are left with two cases:
\begin{itemize}
\item[(A')] The rows of the  diagram obtained by removing the
rightmost square of the $(i+1)$-th row have distinct lengths. In
this case, we fill the rightmost square with $t$ and the remaining
blank squares of the $(i+1)$-th row with $t-1$.

\item[(B')] The removal of the rightmost square of the $(i+1)$-th row
does not result in a bar tableau. Suppose that the $(i+1)$-th row
has $m$ squares in total. It can only happen that the row underneath
the $(i+1)$-th row has $m-1$ squares and all these squares are
filled with $0$'s. By interchanging the location of $0$'s in these
two rows, we get a new configuration $\mathcal{C}^{\prime}$, see
Figure \ref{case2'}. From Lemma \ref{varyzero2-3} we deduce that
$\kappa({\mathcal{C}^*})=\kappa({\mathcal{C}^{\prime}})$. So we can
transform ${\mathcal{C}^*}$ to ${\mathcal{C}'}$ and continue to fill
up the $(i+1)$-th row.
\end{itemize}

Finally, we arrive at  a  shifted diagram  whose rows are all filled
up. Clearly, for those rows containing at least one $0$ there are
$o_s^*+2e_s^*$ bars that are generated in the construction, and for
those rows containing no $0$'s there are
$\max(o_r^*,e_r^*+((e_r^*+o_r^*)\ \mathrm{mod}\ 2))$ bars that are
generated. It has been shown that during the procedure of filling
the diagram with nonnegative numbers if the configuration
${\mathcal{C}^*}$ is transformed to another configuration
${\mathcal{C}^{\prime}}$, then $\kappa({\mathcal{C}^\prime})$
remains equal to $\kappa({\mathcal{C}^*})$. Hence the above
procedure leads to a skew bar tableau of shape $\lambda/\mu$ with
$\kappa({\mathcal{C}^*})$ bars. This completes the proof. \qed

\vspace{.2cm} \noindent{\bf Acknowledgments.} This work was
supported by  the 973 Project, the PCSIRT Project of the Ministry of
Education, the Ministry of Science and Technology, and the National
Science Foundation of China.

\end{document}